\documentclass{article}
\usepackage{graphicx}
\usepackage{amsmath}

\input{tcilatex}

\begin{document}

\title{C\'{a}lculo do Tensor M\'{e}trico de Sistemas Din\^{a}micos Lineares via
Computa\c{c}\~{a}o Alg\'{e}brica}
\author{Jo\~{a}o Jos\'{e} de Farias Neto \\
%EndAName
Instituto de Estudos Avan\c{c}ados - CTA\\
Divis\~{a}o de Inform\'{a}tica\\
Rod. dos Tamoios, km 5,5, S\~{a}o Jos\'{e} dos Campos, SP\\
joaojfn@ieav.cta.br}
\maketitle

As refer\^{e}ncias [1], [2] e [3] descrevem a aplica\c{c}\~{a}o de geometria
diferencial ao problema da identifica\c{c}\~{a}o de sistemas din\^{a}micos
lineares, isto \'{e}, a determina\c{c}\~{a}o dos par\^{a}metros que os
definem a partir de amostras de s\'{e}ries de entrada e de sa\'{i}da. O
conjunto de tais sistemas com o mesmo n\'{u}mero de s\'{e}ries de entrada e
de sa\'{i}da e mesmo grau de McMillan constitui uma variedade diferenci\'{a}%
vel riemaniana. O interesse em determinar-se o tensor m\'{e}trico em cada
ponto de tais variedades tem dois motivos: um, puramente matem\'{a}tico, de
investiga\c{c}\~{a}o das propriedades geom\'{e}tricas daquelas variedades, e
outro, computacional, de melhoria dos algoritmos de otimiza\c{c}\~{a}o
usados no processo de estima\c{c}\~{a}o dos par\^{a}metros de sistemas
reais. Os algoritmos do tipo quase-Newton usam uma aproxima\c{c}\~{a}o da
hessiana da fun\c{c}\~{a}o objetivo (fun\c{c}\~{a}o de verossimilhan\c{c}a,
por exemplo); se essa aproxima\c{c}\~{a}o for substitu\'{i}da pelo tensor m%
\'{e}trico G($\theta )$ correspondente ao ponto $\theta $ em que se est\'{a}%
, a trajet\'{o}ria percorrida at\'{e} chegar-se ao ponto \'{o}timo fica
quase independente da parametriza\c{c}\~{a}o usada para descrever o sistema;
consegue-se, assim, fazer \ \ \textit{identifica\c{c}\~{a}o independente de
coordenadas.}

Seja H(z) a fun\c{c}\~{a}o de transfer\^{e}ncia de um sistema linear a tempo
discreto com m entradas, m sa\'{i}das e grau de McMillan n. H(z) \'{e} uma
matriz de fun\c{c}\~{o}es racionais de z. Suponha-se o sistema est\'{a}vel.
Qualquer sistema desse tipo pode ser descrito por 2mn par\^{a}metros $\theta
_{i}$ , i=1,2...,2mn, de modo que H(z) \'{e} uma fun\c{c}\~{a}o do vetor $%
\theta $ . Em [5], o autor deste artigo prova que o elemento g$_{ij}$ do
tensor pode ser obtido pela f\'{o}rmula:

\begin{equation*}
\underline{\fbox{g$_{ij}=\frac{1}{2\pi \text{i}}\doint_{C}tr\left[ \frac{%
\partial H(z)}{\partial \theta _{i}}\frac{\partial H^{T}(z^{-1})}{\partial
\theta _{j}}\right] z^{-1}dz$}}
\end{equation*}

onde i=$\sqrt{-1}$ , C \'{e} o c\'{i}rculo unit\'{a}rio centrado na origem
do plano complexo e tr representa tra\c{c}o.

Nas parametriza\c{c}\~{o}es superpostas ou pseudo-can\^{o}nicas, as
componentes de $\theta $ s\~{a}o naturalmente divididas em dois grupos, cada
um com mn componentes: as que s\'{o} aparecem no lado esquerdo de um modelo
ARMA (ou, equivalentemente, na matriz A que multiplica o estado em um modelo
em espa\c{c}o de estados) e as outras. Seja I o primeiro grupo e J, o
segundo.

No modelo ARMA(p,p)

$%
A_{0}y_{t}+A_{1}y_{t-1}+...+A_{p}y_{t-p}=B_{0}u_{t}+B_{1}u_{t-1}+...+B_{p}u_{t-p}
$

onde y$_{t}$,u$_{t}$ $\in R^{m}$ \ e \ A$_{i}$ e B$_{i}$ s\~{a}o matrizes $%
m\times m$ , tem-se

H(z)=A$^{-1}\left( z\right) B(z)$

onde \ \ \ \ \ \ \ \ \ A(z)= A$_{0}z^{p}+A_{1}z^{p-1}+...+A_{p}\mathstrut
\strut $

e \ \ \ \ \ \ \ \ \ \ \ \ \ \ \ B(z)= B$_{0}z^{p}+B_{1}z^{p-1}+...+B_{p}\ $

Sejam M(z) e K as matrizes da parametriza\c{c}\~{a}o superposta ARMA
definida em [4], onde p=max\{n$_{i}$\}, sendo n$_{i}$ os \'{i}ndices de
Kronecker do sistema. Ent\~{a}o, o autor do presente artigo demonstra em [6]
que

\textbf{Caso 1: \ }\ \ \ \ i,j$\in $J

g$_{ij}=\frac{1}{2\pi \text{i}}\doint_{C}tr\left[ A^{-1}(z)M(z)\frac{%
\partial K}{\partial \theta _{i}}\frac{\partial K^{T}}{\partial \theta _{j}}%
M^{T}(z^{-1})A^{-T}(z^{-1})\right] z^{-1}dz$

\textbf{Caso 2: }\ \ \ \ \ i,j$\in $I

g$_{ij}=\frac{1}{2\pi \text{i}}\doint_{C}tr[A^{-1}(z)\left( \frac{\partial
M(z)}{\partial \theta _{i}}K+\frac{\partial A(z)}{\partial \theta _{i}}\left[
I-H(z)\right] \right) .$

.$\left( K^{T}\frac{\partial M^{T}(z^{-1})}{\partial \theta _{j}}+\left[
I-H^{T}(z^{-1})\right] \frac{\partial A^{T}(z^{-1})}{\partial \theta _{j}}%
\right) A^{-T}(z^{-1})]z^{-1}dz$

\textbf{Caso 3: \ }\ \ \ \ i$\in $I \ \ and \ \ j$\in $J \ \ \ 

g$_{ij}=\frac{1}{2\pi \text{i}}\doint_{C}tr$

$\left[ A^{-1}(z)\left( \frac{\partial M(z)}{\partial \theta _{i}}K+\frac{%
\partial A(z)}{\partial \theta _{i}}\left[ I-H(z)\right] \right) \frac{%
\partial K^{T}}{\partial \theta _{j}}M^{T}(z^{-1})A^{-T}(z^{-1})\right]
z^{-1}dz$

No modelo em espa\c{c}o de estados

$x_{t+1}=Ax_{t}+Bu_{t}$

$y_{t}=Cx_{t}+u_{t}$

onde $x_{t}\in R^{n}$ and $u_{t},y_{t}\in R^{m}$

com A, B e C dadas em fun\c{c}\~{a}o de $\theta $ pela parametriza\c{c}%
\~{a}o superposta definida em [4], demonstra-se em [6] que:

\textbf{Caso 1}: \ \ \ \ i,j$\in $J \ \ 

$g_{ij}=\frac{1}{2\pi \text{i}}\doint_{C}tr\left[ C(zI-A)^{-1}\frac{\partial
B}{\partial \theta _{i}}\frac{\partial B^{T}}{\partial \theta j}%
(z^{-1}I-A)^{-T}C^{T}\right] z^{-1}dz$

\textbf{Caso 2: \ }\ i,j$\in $I

$g_{ij}=\frac{1}{2\pi \text{i}}\doint_{C}tr$

$\left[ C(zI-A)^{-1}\frac{\partial A}{\partial \theta _{i}}%
(zI-A)^{-1}(z^{-1}I-A)^{-T}\frac{\partial A^{T}}{\partial \theta _{j}}%
(z^{-1}I-A)^{-T}C^{T}\right] z^{-1}dz$

\textbf{Caso 3:} \ \ \ i$\in $I \ \ and \ \ j$\in $J \ \ 

$g_{ij}=\frac{1}{2\pi \text{i}}\doint_{C}tr\left[ C(zI-A)^{-1}\frac{\partial
A}{\partial \theta _{i}}(zI-A)^{-1}\frac{\partial B^{T}}{\partial \theta j}%
(z^{-1}I-A)^{-T}C^{T}\right] z^{-1}dz$

No caso de sistemas estoc\'{a}sticos, isto \'{e}, com entradas ruidosas
desconhecidas (o que torna a sa\'{i}da probabil\'{i}stica), sejam $\Gamma
_{i}=E(y_{t}y_{t+i}^{T}),i=0,1,2,...$. Ent\~{a}o, este autor mostra, em [7],
que as seguintes duas f\'{o}rmulas fornecem tensores, conforme a defini\c{c}%
\~{a}o que se escolha para o produto interno no espa\c{c}o tangente \`{a}
variedade:

$g_{ij}=\frac{1}{2\pi \text{i}}\doint_{C}tr\left[ \frac{\partial U(z)}{%
\partial \theta _{i}}\frac{\partial U^{T}(z^{-1})}{\partial \theta _{j}}%
\right] z^{-1}dz$

onde:$U(z)=\sum_{i=0}^{\infty }\Gamma _{i}z^{-i}$

Ou:

$g_{ij}=\frac{1}{2\pi \text{i}}\doint_{C}tr\left[ \frac{\partial T(z)}{%
\partial \theta _{i}}\frac{\partial T^{T}(z^{-1})}{\partial \theta _{j}}%
\right] z^{-1}dz$

onde T(z) = $\sum_{i=-\infty }^{\infty }\Gamma _{i}z^{-i}$ = H(z)RH$^{T}$(z$%
^{-1}$) \'{e} a fun\c{c}\~{a}o densidade espectral da sa\'{i}da $\left\{
y_{t}\right\} $ do sistema din\^{a}mico linear cuja fun\c{c}\~{a}o de
transfer\^{e}ncia \'{e} H(z) e cuja entrada \'{e} um ru\'{i}do branco $%
\left\{ \varepsilon _{t}\right\} $, de m\'{e}dia nula e covari\^{a}ncia R.

A conclus\~{a}o \'{e} que utilizando-se uma linguagem de programa\c{c}\~{a}o
que contenha comandos de computa\c{c}\~{a}o alg\'{e}brica, como MAPLE ou
MATHEMATICA, pode-se, usando as f\'{o}rmulas aqui apresentadas, calcular o
tensor e incorpor\'{a}-lo ao algoritmo de estima\c{c}\~{a}o dos
par\^{a}metros ou simplesmente us\'{a}-lo para descrever as propriedades
geom\'{e}tricas de conjuntos de sistemas lineares din\^{a}micos. Embora as
f\'{o}rmulas resultantes cres\c{c}am rapidamente em complexidade, a
hip\'{o}tese de estabilidade permite que muitos termos de ordem superior
sejam desprezados, obtendo-se aproxima\c{c}\~{o}es ainda \'{u}teis.

\textbf{Refer\^{e}ncias:}

[1]B. Hanzon. Identifiability, recursive identification and spaces of linear
dynamical systems. Tese de doutorado, Departamento de Econometria,
Universidade Erasmus, Rotterdam,1986.

[2] R. Peeters. System identification based on Riemannian geometry: theory
and algorithms.Tese de doutorado. Free University of Amsterdam, 1994, e
Research Report nr. 64, Tinbergen Institute Research Series, Tinbergen
Institute, Rotterdam.

[3]Chou, T. C., Geometry of Linear Systems and Identification. Tese de
doutorado.Trinity College, Cambridge, 1994.

[4] M. Gevers, V. Wertz. Uniquely identifiable state-space and ARMA
parametrizations for multivariable systems. Automatica, vol. 20, n. 3, 1984,
333-347.

[5] Farias Neto, J.J. Tensor m\'{e}trico riemaniano ARMA. Relat\'{o}rio de
Pesquisa. Departamento de Matem\'{a}tica, Universidade do Estado de Santa
Catarina.

http://www.geocities.com/joaojfn. 1998.

[6] Farias Neto, J.J.. A note on the computation of linear systems metric
tensors. Relat\'{o}rio de Pesquisa, Departamento de Matem\'{a}tica,
Universidade do Estado de Santa Catarina. http://www.geocities.com/joaojfn.
1999.

[7] Farias Neto, J.J. Tensor m\'{e}trico riemaniano para sistemas estoc\'{a}%
sticos. Relat\'{o}rio de Pesquisa, Departamento de Matem\'{a}tica,
Universidade do Estado de Santa Catarina. http://www.geocities.com/joaojfn.
2000.

\end{document}